\documentclass{article}
\usepackage{rukshin}

\begin{document}
\title{Уникальносоставленные фигуры,}
\author{А.~М.~Петрунин, С.~Е.~Рукшин
}
\date{}
\maketitle

По отношению к равносоставленности некоторые выпуклые фигуры обладают одним любопытным свойством: например,  если круг $K$ равносоставлен выпуклой фигурe $F$, то $F$ конгруэнтна $K$.
Дубинс, Каруш и Хирш доказали, что тем же свойством обладает и эллипс (см.  \cite{DKH}). 
Мы же даём полную классификацию таких фигур. 

Эти результаты были получены  и даже написаны нами почти 20 лет назад, тогда первый из авторов занимался в математическом кружке под руководством второго. 
Вскоре после написания статья (напечатанная на печатной машинке в единственном экземпляре) была утеряна, 
даже не дойдя до редакции журнала.

Мы попытались переписать её как можно понятнее. 
Для упрощения мы доказывем последовательно несколько теорем, каждая следующая из которых обобщает предыдущую.
Вам потребуется только знание доказательсва теоремы Бойяи --- Гервина, которое можно найти, например, в  книжках Болтянского \cite{B1} или \cite{B2}; 
в пятой части мы также используем некоторые результаты вещественного анализа (достаточно первых четырёх глав книжки Натансона \cite{N}).

\section{Обозначения и определения.}

\begin{itemize}
\item $F$ будет обозначать плоскую и чаще всего выпуклую фигуру,
\item $S(F)$  обозначает площадь $F$,
\item $\partial F$ обозначает  кривую, ограничивающую $F$
\end{itemize}

\begin{thm}{Определение}
 Две плоских фигуры $F$ и $G$ называются \emph{равносоставленными} (или $F\sim G$), если одну можно разрезать отрезками прямой на конечное число кусков и составить из этих кусков другую фигуру.
\end{thm}

\begin{thm}{Определение}
Выпуклая фигура $F$ называется \emph{уникальносоставленной}, если любая выпуклая фигура, равносоставленная $F$, является конгруэнтной $F$. 
\end{thm}

\begin{thm}{Определение}
Две кривые (или два набора кривых) $\alpha$ и  $\beta$ называются  \emph{равносоставленными} (или $\alpha\sim\beta$), если первую(ый) можно разбить на конечное число дуг и составить из них вторую(ой). 
\end{thm}

\begin{thm}{Определение}
Две кривые $\alpha$ и  $\beta$ называются \emph{стабильно равносоставлеными} (или $\alpha\approx\beta$), если  из них можно исключить конечное число  отрезков прямой так, что остaвшиеся два набора кривых  равносоставлены.
\end{thm}

Следующее утверждение является простым обобщением теоремы Бойяи --- Гервина утверждающей что равновеликие многоугольники равносоствалены, (см. \cite{B1} или \cite{B2}):

\begin{thm}{Утверждение}\label{cl-1}
Две выпуклые фигуры $F$ и $G$ равносоставлены тогда и только тогда, когда  их площади равны и кривые, их ограничивающие, стабильно равносоставлены. Или
$$F\sim G\ \Leftrightarrow\ S(F) = S(G)\  \mbox{и}\ \partial F\approx\partial G.$$

\end{thm}

\parit{Доказательство; необходимость:}  совокупность всех кривых, ограничивающих куски разбиения, состоит из дуг первой кривой плюс какое-то количество «внутренних» отрезков прямой, и из тех же кусков состоит и вторая фигура. Таким образом,  если из границы $F$ вырезать отрезки прямой, которые  соответствуют «внутренним» отрезкам разбиения $G$, и аналоагично поступить с границей $G$, то оставшиеся куски границ будут равносоставлены.

\parit{Достаточность:} рассмотрим разбиение границ $F$ и $G$ как в  определении стабильной равносоставленности. 
Если к каждой такой дуге фигур $F$ и $G$ провести хорду, то образовавшиеся горбушки (сегменты) для равных дуг равны.
Так как $S(F)=S(G)$ то и площади оставшихся многоугольников равны, а значит, они равносоставлены по теореме Бойяи --- Гервина.
\qeds

\section{Уникальносоставленность круга.}

\begin{thm}{Теорема}\label{thm-1}
Круг уникальносоставлен.
\end{thm}

Для доказательства нам потребуется  следующая лемма:

\begin{thm}{Лемма}\label{lm-1}
Если граница  выпуклой фигуры $F$ состоит из конечного числа отрезков прямой и дуг окружностей радиуса $R$ с общей угловой мерой $360^\circ$, то она является $R$-окрестностью выпуклого многоугольника. 

В частности,  $F$ содержит в себе круг радиуса $R$.
\end{thm}

\parit{Доказательство теоремы \ref{thm-1} по модулю леммы \ref{lm-1}.}
Обозначим через $K$ круг радиуса  $R$. Пусть $F$ есть выпуклая фигура, такая что $K\sim F$. Тогда из утверждения~\ref{cl-1} мы имеем: $S(F)=S(K)$  и  граница $F$ {\it стабильно равносоставлена} окружности радиуса $R$. 
В частности, $F$ удовлетворяет лемме \ref{lm-1}. Таким образом,
$F$ содержит в себе копию $K$. 
Так как $S(F)=S(K)$,  получаем $K\cong F$.

\rightline{$\square$}

В доказательстве леммы мы будем пользоваться следующим фактом, который предоставляем читателю в качестве упражнения:

\begin{thm}{Упражнение}\label{ex-1}
Пусть $A_1A_2\cdots A_n$ есть замкнутая ломаная, которая при обходе по ней поворачивает всё время против часовой стрелки, и общий угол поворота равен $360^\circ$ --- тогда  ломаная $A_1A_2\cdots A_n$ ограничивает выпуклый многоугольник.
\end{thm}

\parit{Доказательство леммы \ref{lm-1}.}
Пусть $\alpha$ есть выпуклая кривая стабильно равносоставленная окружности радиуса $R$.
Заметим, что в этом случае кривая не имеет угловых точек. 
Действительно,  при обходе вокруг кривой необходимо повернуть на $360^\circ $, и ровно на этот же угол мы поворачиваем, пройдя по всем дугам окружности; таким образом, на углы не остаётся места (то есть, если бы имелась ещё и угловая точка, то общий угол поворота был бы больше чем $360^\circ $, что невозможно).
 
\begin{wrapfigure}{l}{40mm}
\begin{lpic}[t(0mm),b(0mm),r(0mm),l(0mm)]{R-nbhd-kr-zel-zh(0.30)}
\lbl{62,73;$M$}
\end{lpic}
\end{wrapfigure}

Таким образом, кривая $\alpha$ состоит из поочерёдно сменяющихся дуг окружностей и отрезков прямых, так что отрезки и дуги продолжают друг друга в том же направлении.
Для каждого отрезка прямой в $\alpha$ рассмотрим его параллельный перенос в перпендикулярном направлении внутрь фигуры на расстояние $R$.
При этом концы отрезков, соседниe через дугу, перейдут в одну точку, и мы получим замкнутую ломаную из параллельных переносов всех отрезков $\alpha$.
Звенья этой ломаной поворачивают в одну и ту же сторону (так же как и $\alpha$), и, так же как у $\alpha$, общий угол поворота будет $360^\circ$. 

Воспользовавшись упражнением~\ref{ex-1}, мы получаем, что эта ломаная ограничивает выпуклый многоугольник. Oбозначим его $M$, и из построения легко видеть, что сама фигура $F$ есть множество точек на расстоянии $R$ от $M$.

\section{Уникальносоставленность линзы.}

\begin{thm}{Теорема}\label{thm-2}
Пересечение двух кругов одинакового радиуса есть уникальносоставленная фигура.
\end{thm}

Для краткости, давайте назовём пересечение двух кругов радиуса $R$ \emph{линзой} и будем обозначать её $L_\omega$, где $\omega$ есть угловая мера её дуг (в частности, $L_{2{\cdot}\pi}$ есть круг радиуса $R$). 
С помощью точно такого же рассуждения, что и в доказательстве теоремы \ref{thm-1}, теорема \ref{thm-2} сводится к  следующей лемме:

\begin{thm}{Лемма}\label{lm-2}
 Если выпуклая фигура $F$ имеет границу, стабильно равносоставленную границе линзы
  $L_\omega$,  то $S(F)\geqslant S(L_\omega)$. Более того, если $S(F)=S(L_\omega)$, то  $F\cong L_\omega$.
\end{thm}

Доказательство будет проведено в два шага. Сначала мы докажем следующее, более слабое утверждение, а уже потом приступим к общему случаю.

\begin{thm}{Лемма}\label{lm-3}
  Если центрально-симметричная выпуклая фигура $F$ имеет границу, стабильно равносоставленую границе линзы $L_\omega$,  то $S(F)\geqslant S(L_\omega)$.
Более того, если $S(F)=S(L_\omega)$, то  $F\cong L_\omega$.
 \end{thm}

\parit{Доказательство леммы \ref{lm-3}.} 
В доказательстве мы воспользуемся двумя процедурами: «вырезанием параллелограмма» и «четырёхшарнирным сдвигом». 
Первый (применённый несколько раз) позволит убрать из границы $F$ все отрезки прямых.
Второй (также применённый несколько раз) позволит убрать все угловые точки $F$ кроме двух.

\begin{wrapfigure}{r}{54mm}
\begin{lpic}[t(-5mm),b(0mm),r(0mm),l(0mm)]{parallelogramm-kr(0.33)}
\lbl[l]{73,50;$A'$}
\lbl[r]{5,24;$A$}
\lbl[bl]{65,70;$B'$}
\lbl[r]{11,4;$B$}
\end{lpic}
\end{wrapfigure}

\parit{Вырезаниe параллелограмма.}
Предположим, на границе $F$ есть отрезок прямой $AB$, и пусть $A'B'$ обозначает центрально-симметричный ему отрезок. 
Тогда из $F$ можно вырезать параллелограмм $ABA'B'$ и из 
оставшихся двух кусков составить центрально-симметричную выпуклую фигуру $F'$. При этом,  очевидно, $S(F')<S(F)$ и $\partial F'\approx \partial F$ (т.е. их границы стабильно равносоставлены).

Повторив эту операцию столько раз, сколько возможно, выбирая каждый раз новую пару
отрезков, мы получаем  центрально-симметричную выпуклую фигуру
  $F_1$ без отрезков прямой на границе такую, что
$S(F_1)\le S(F)\ \ \mbox{и}\ \ \partial F_1\approx \partial F.$

\begin{wrapfigure}{r}{35mm} 
\begin{lpic}[t(-5mm),b(0mm),r(0mm),l(0mm)]{chet-shar-kr(0.35)}
\lbl[b]{43,155;$B$}
\lbl[t]{50,82;$B'$}
\lbl[r]{4,121;$A$}
\lbl[l]{88,118;$A'$}
\end{lpic}
\end{wrapfigure}

\parit{Четырёхшарнирный сдвиг.} 
Пусть на границе $F$ лежат две пары центрально-симметричных угловых точек $A,A' $ и $B,B' $. 
Разрезав  $F$ отрезками $AB,\ BA',\ A'B'$ и $B'A $, мы получаем параллелограм $ABA'B'$ и четыре горбушки (см. рисунок).
Не умоляя общности, можно считать, что угол $\angle ABA'$ тупой. 
Пусть $\theta$ есть внешний угол границы $F$ при $B$. 
Представьте, что в вершины параллелограмма вставлены шарниры, а стороны сделаны из жёсткого материала, при этом горбушки жёстко приделаны к сторонам. 
Давайте повернём сторону $AB$ вокруг $B$ на угол $\theta$  так,  чтобы  угол  $\angle ABA'$   увеличился. 
Тогда площадь  параллелограмма $ABA'B'$, а значит и площадь фигуры, образованной параллелограммом и четырьмя горбушками, уменьшится; при этом точки $B$ и $B'$ перестанут быть угловыми. 
Обозначим полученную фигуру $F'$

Начнём с фигуры $F_1$ и проделаем с ней четырёхшарнирный сдвиг столько раз, сколько возможно, выбирая каждый раз новую пару пар центрально-симметричных угловых точек, мы получим фигуру $F_2$ такую, что:

а) $F_2$ имеет не более одной пары угловых точек,

б) $F_2$ центрально-симметрична,

в)  $S(F_2)\le S(F_1)$.

д) Граница $F_2$ равносоставлена границе  $F_1$ и у $F_2$ нет отрезков прямой на границе.

Из этих условий нетрудно видеть, что $F_2\cong L_\omega$. 
Таким образом,
$$S(F)\ge S(F_1)\ge S(F_2)=S(L_\omega).$$

\ 

Отметим, что $F_2$ была получена из $F$ как результат последовательности шагов, в каждом из которых площадь уменьшалась. 
Значит, в случае равенства $S(F)=S(L_\omega)$ мы вовсе не делали ни одного шага, то есть $F\cong F_1\cong F_2\cong L_\omega$.\qeds

\parit{Доказательство леммы \ref{lm-2}.}
Пусть $F$ есть фигура описанного типа. 
Проведём две параллельные прямые $a$ и $b$ так, что вся $F$ лежит в полосе между ними, и обе прямые касаются $F$.
Обозначим через $A$ и $B$  точки пересечения  $F$  соответственно с $a$ и $b$ (если прямая соприкасается с $F$ вдоль отрезка, то следует выбрать произвольную точку на отрезке). 
Пусть $\alpha$ и $\beta$ обозначают соответственно общую угловую величину всех дуг окружностей на участках границы $F$ от $A$ до $B$ и от $B$ до $A$ против часовой стрeлки.

\begin{thm}{Упражнение}
Покажите,  что пару $A,B$ можно выбрать с дополнительным свойством что  $\alpha=\beta$. 
\end{thm}

\parit{Подсказка.} Для этого надо просто прокрутить пару $A,B$ непрерывно вокруг  $F$, так что $A$ перешло бы в $B$, а $B$ в $A$, и воспользоваться тем, что непрерывная функция на отрезке принимает все промежуточные значения. 
Мы также советуем посмотреть статью Болтянского и Савина \cite{BS}, в которой обсуждаются геометрические задачи, в решении которых применяется эта идея.

\vskip0.2cm

Теперь мы можем считать, что $\alpha=\beta=\omega/2$.
Разрежем $F$ отрезком $AB$ на две части $F_-$ и $F_+$. Из двух копий $F_-$ можно составить центрально-симметричную фигуру $\widetilde F_-$.
Аналогично поступим с $F_+$ и получим фигуру
$\widetilde F_+$.

Согласно лемме \ref{lm-3} каждая из фигур $\widetilde F_-$ и $\widetilde F_+$ имеет площадь не меньше, чем линза $L_{\omega}$. 
А значит 
$$S(F)=S(F_+)+S(F_-)=\tfrac{1}{2}{\cdot}\!\left(S(\widetilde F_+)+S(\widetilde F_-)\right)\geqslant S(L_{\omega}).$$

\begin{center}
\begin{lpic}[t(0mm),b(0mm),r(0mm),l(0mm)]{strong2-kr-zel(0.35)}
\end{lpic}
\end{center}

Отметим, что равенство может достигаться только в случае, если обе фигуры 
$\widetilde F_-$ и $\widetilde F_+$ конгруэнтны $L_{\omega}$, а значит $F$ имеет не более двух угловых точек и не имеет отрезков прямой на границе, а значит $F\cong L_{\omega}$.

\rightline{$\square$}

\section{Круговые многоугольники}

Мы приступаем к случаю круговых многоугольников, начнём с определений:

\begin{thm}{Определение}
Выпуклая фигура называется круговым многоугольником, если её граница равносоставлена конечному набору дуг окружностей и отрезков прямой.
\end{thm}

\begin{thm}{Определение}\label{oval}
Круговой многоугольник $M$ называется овальным (или овалом), если он обладает следующими свойствами:

а) не имеет отрезков прямой на границе;

\begin{wrapfigure}{r}{50mm} 
\includegraphics[scale=0.35]{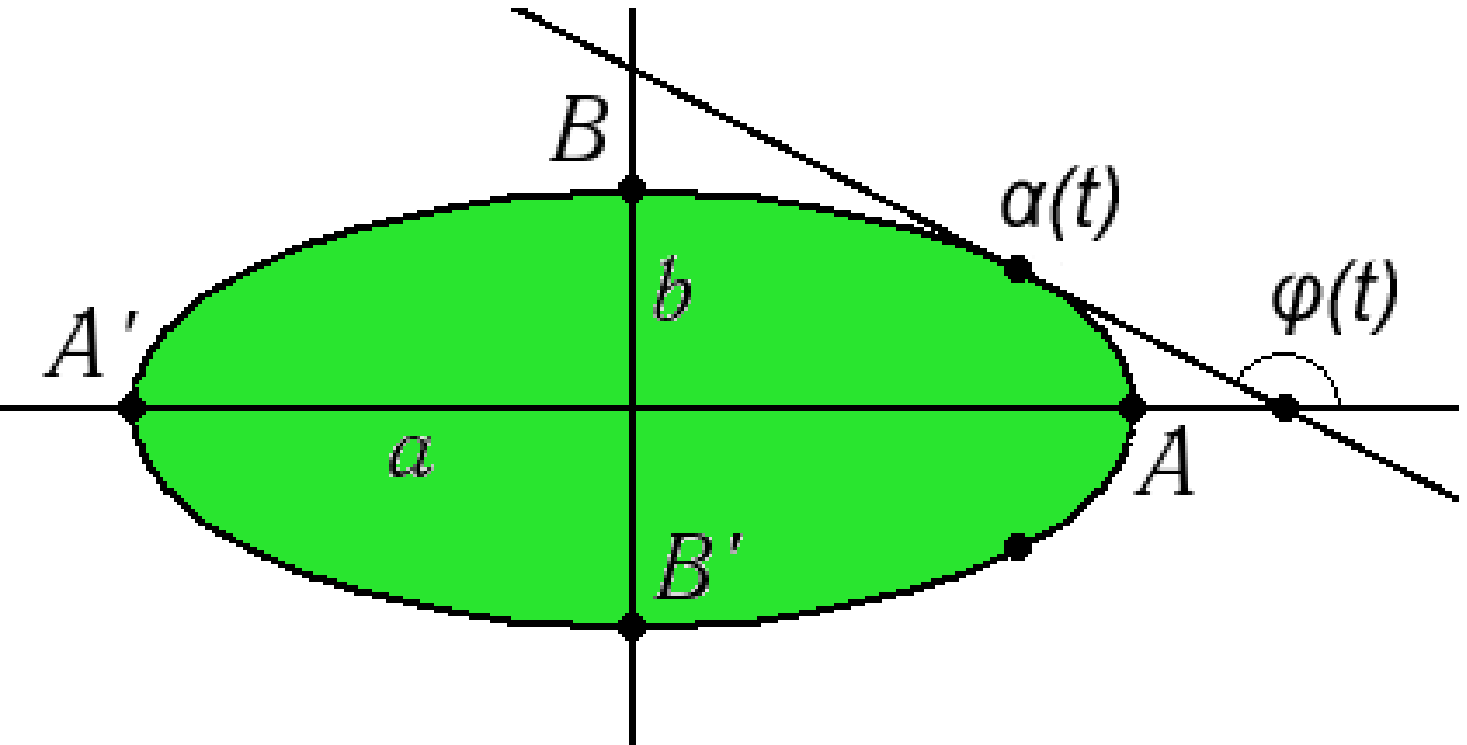}
\end{wrapfigure}

б) имеет две перпендикулярные оси симметрии $a$ и $b$. Давайте обозначим через $A, A'$ и $B,B'$ их точки пересечения с границей $M$;

в) если граница $M$ имеет угловые точки то они находятся на оси $a$, т.е. в точках $A$ и $A'$;

г) радиусы дуг окружностей, ограничивающей  $M$, монотонно возрастают  от $A$ до $B$. 
\end{thm}

Легко видеть, что  круговой многоугольник овален тогда и только тогда, когда овальна его круговая окрестность.

\begin{thm}{Теорема}\label{thm-3}
Круговой многоугольник  уникальносоставлен тогда и только тогда, когда он овален.
\end{thm}

В следующем разделе мы обобщим этот результат на случай всех выпуклых фигур.

Нетрудно видеть, что в классе круговых многоугольников со стабильно равносоставленной границей есть единственный, с точностью до конгруэнтности, овал. 
Как и в случаях круга и линзы, теорема \ref{lm-3} для круговых многоугольников немедленно следует из аналога леммы \ref{lm-2}:

\begin{thm}{Лемма} \label{lm-4}
Если круговой многоугольник $M$ имеет границу, стабильно равносоставленую границе овала $V$,  то $$S(M)\geqslant S(V).$$
Более того, если $S(M)=S(V)$, то  $M\cong V$.
\end{thm}

Как и в случае с линзой, доказательство будет проведено в два шага; сначала мы докажем аналогичное более слабое утверждение:

\begin{thm}{Лемма} \label{lm-5}
Если центрально-симметричный круговой многоугольник $M$ имеет границу, стабильно равносоставленую границе овала $V$,  то 
$$S(M)\geqslant S(V).$$
Более того, если $S(M)=S(V)$ то  $M\cong V$.
\end{thm}

Доказательство леммы \ref{lm-5} проходит индукцией по количеству дуг различных радиусов на границе $M$, при этом мы считаем угловые точки за дуги нулевого радиуса. 
Доказательство для круга можно принять за базу индукции, а доказательство для линзы можно рассматривать как часть первого шага в этой индукции (на границе линзы присутствуют дуги радиуса $R$ и радиуса $0$). Тем не менее, в доказательстве участвует ещё пара дополнительных идей, которые нам не были нужны раньше.

В частности, нам потребуются следующие утверждения, доказательство которых мы предоставляем читателю:

\begin{thm}{Упражнение} \label{ex-nbhd}
Пусть $F$ и $G$ суть выпуклые фигуры с равносоставленной границей, a $F_R$ и $G_R$ --- их $R$-окрестности; тогда $F \sim G$ тогда и только тогда, когда $F_R \sim G_R$.
\end{thm}

\begin{thm}{Упражнение}\label{max}
Если выпуклая фигура $F$ имеет границу, равносоставленую границе овала $V$, 
то диаметр $F$ не превосходит диаметра $V$. 
Более того, в случае равенства диаметров, $F$ конгруэнтна $V$.
\end{thm}

\parit{Доказательство леммы \ref{lm-5}.} Пусть $M$ есть круговой многоугольник с $n$ 
различными радиусами дуг на границе. 

Как и в случае линзы, применив несколько раз «вырезаниe параллелограмма», мы сводим задачу к случаю, когда на границе $M$ отсутствуют отрезки прямой. В частности, мы можем считать что граница $M$, уже не только стабильно равносоставлена, но и равносоставлена границе $V$.

Если на границе кругового  многоугольника $M$ нет угловых точек, то существует круговой  многоугольник $M'$, такой что
$M=M'_R$ (напомним, что $M'_R$ обозначает $R$-окрестность $M'$). 
Таким образом, из упражнения \ref{ex-nbhd} следует что достаточно доказать лемму для случая, когда на границе $M$ есть угловые точки.

Eсли у $M$ более одной пары угловых точек, то можно применить несколько раз «четырёхшарнирный метод» и свести задачу к случаю, когда на границе $M$ есть ровно две угловые точки.

Далее, пусть $A,\, A'$ есть едиственная пара угловых точек $M$, 
пусть $\alpha$ есть внешний угол при каждой из этих точек, и $r_1$ есть минимальный ненулевой радиус дуг на границе $M$. 

\begin{wrapfigure}{r}{30mm}
\begin{lpic}[t(-5mm),b(0mm),r(0mm),l(0mm)]{weak-kr-zel(0.40)}
\end{lpic}
\end{wrapfigure}

Разрежeм границу $M$ на две дуги по точкам $A$ и $ A'$ и продолжим оба 
 конца первой из дуг двумя дугами радиуса $r_1$ с угловой величиной $\alpha$. К полученной дуге можно приставить вторую из дуг от  $A$ до $ A'$ до образования выпуклой замкнутой кривой без угловых точек. Обозначим через $\bar M$ фигуру, ограниченную полученной кривой.

Проделаем ту же операцию с овалом $V$, получим овал $\bar V$ с границей равносоставленной $\bar M$. 

Заметим что при этом 
$$S(\bar V)-S(V)\ge S(\bar M)-S(M),$$
причём равенство достигается только если $M$ конгруэнтна $V$.
Действительно, $V$ получается из $\bar V$ добавлением двух круговых сегментов и прямоугольника с одной стороной, равной диаметру $\bar V$, и другой, равной $2r_1\sin(\alpha/2)$. 
С другой стороны, $M$ получается из $\bar M$ добавлением тех же двух сегментов и параллелограмма с одной стороной, не превосходящей диаметра $\bar V$ (в силу упражнения \ref{max}), и другой стороной, также равной $2r_1\sin(\alpha/2)$.

Теперь мы получили круговой многоугольник $\bar M$ с $n-1$ различными радиусами дуг на границе. 
По предположению индукции, 
$$S(\bar M)\ge S(\bar V).$$ 
Сложив эти неравенства, получаем
$$S(M)\ge S(V).$$

Случай равенства возможен, только если диаметер $M$ равен диаметру $V$.   
Значит, в силу упражнения \ref{max}, $M\cong V.$\qed

\medskip

Следующее доказательство
схоже на  доказательство леммы \ref{lm-2} для 
случая линзы, но в нём участвует ещё один забавный трюк.

\parit{Доказательство леммы~\ref{lm-4}.} Пусть $M$ есть круговой многоугольник. Пусть $A$ и $B$  суть точки на границе $M$, такие, что через них можно провести параллельные прямые так, что вся $M$  будет лежать в полосе между ними. Разрежем $M$ отрезком $AB$ на две части $M_-$ и $M_+$. 
Из двух копий $M_-$ составим центрально-симметричную выпуклую фигуру $\widetilde M_-$.
Аналогично поступим с $M_+$ и получим фигуру $\widetilde M_+$. 
Давайте обозначим через $\widetilde V_+$ и $\widetilde V_-$ овалы с границей, стабильно равносоставленной соответственно $\widetilde M_+$ и $\widetilde M_-$

Диаметры овалов $\widetilde V_+$ и $\widetilde V_-$ непрерывно зависят от выбора точек $A$ и $B$. 
Таким образом, прокручиванием $A$ и $B$ вокруг $M$, также как в доказательстве для линзы, 
можно добиться такого выбора $A$ и $B$, что диаметры $\widetilde V_-$ и $\widetilde V_+$ совпадут.

Разрезав $\widetilde V_-$ двумя осями симметрии, мы получим четыре четвертинки конгруэнтные, скажем  $\widetilde V_{-1/4}$. Аналогично поступим с $\widetilde V_-$, получим четыре четвертинки конгруэнтные  $\widetilde V_{+1/4}$.

Из двух копий $\widetilde V_{-1/4}$ и двух копий $\widetilde V_{+1/4}$ можно составить центрально-симметричный круговой многоугольник $N$, как показано на картинке:
\begin{center}
\includegraphics[scale=0.35]{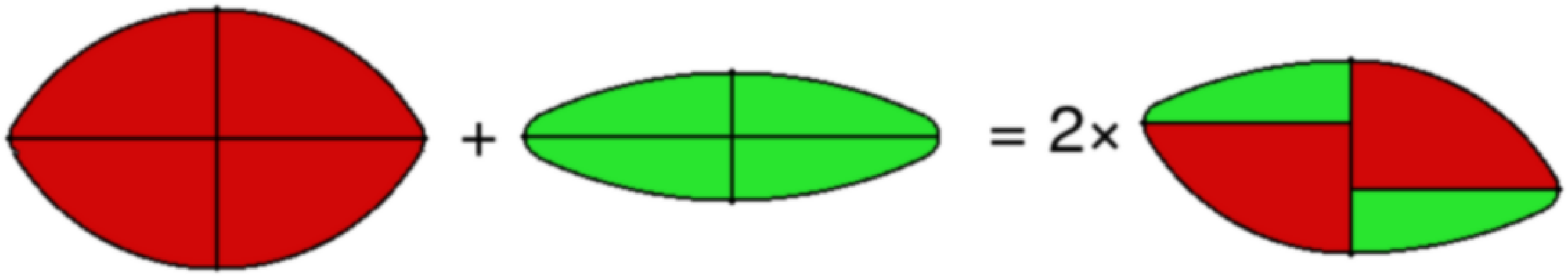}
\end{center}
Применив лемму \ref{lm-5} к $N$ получаем:

$$S(M)=S(M_+)+ S( M_-)=\tfrac12{\cdot}\left(S(\widetilde M_+)+ S(\widetilde M_-)\right)\ge$$
$$\ge \tfrac12{\cdot}\left(S(\widetilde V_+)+ S(\widetilde V_-)\right)=S(N)\ge S(V)$$

Остаётся проверить, что равенство в обоих этих неравенствах влечёт конгруэнтность $M$ и $V$, и это мы щедро предоставляем читателю.

Мы доказали, что любой овальный круговой многоугольник является уникальносоставленной фигурой.  Остаётся ещё проверить,
что не существует других уникальносоставленных  круговых многоугольников.
Применив утверждение~\ref{cl-1}, мы получаем, что
каждый круговой многоугольник равносостaвлен набору из овала и квадрата.  
Таким образом, достаточно предъявить  две неконгруэнтные выпуклые фигуры, каждая из которых равносоставлена нашему набору из квадрата и овала. 
Проверку этого утверждения мы также предоставляем читателю. \qed

\section{Набросок доказательства общего случая.}

Это наиболее технический и наименее геометрический раздел статьи. Для его понимания вам потребуруется знать в некотором объёме теорию функций вещественного переменного (достаточно первых четырёх глав книжки Натансона \cite{N})  и также знать элементарные свойства метрики Хаусдофа, которую мы обозначаем $d_H$, для компактных подмножеств плоскости (см. например, \cite{S}).

Пусть $\alpha:[0,T)\to \RR^2$ есть граница выпуклой фигуры параметризованная длиной в направлении против часовой стрелки. 
У $\alpha$ определена правая производная $\alpha^+(t)$. 
Определим угол поворота кривой $\phi(t)$ как угол от $\alpha^+(0)$ до $\alpha^+(t)$ (измеряемый от $0$ до $2{\cdot}\pi$ против часовой стрелки). 
Из выпуклости $\alpha$ следует, что $\phi(t)$ монотонно возрастающая функция.
Это даёт возможность определить «верхнюю кривизну» $k(t)$ как верхний предел
$$ k(t_0)=\overline{\lim_{t\to t_0}}\,\frac{\phi(t)-\phi(t_0)}{t-t_0}.$$
Так как возрастающая функция имеет производную в почти всех точках, мы получаем, что $k(t)$ конечна для почти всех $t$.

Определим «нижний радиус кривизны» через $R(t)=\frac1{k(t)}$ (считая $0=1/\infty$). 
Это даёт возможность обобщить определение овала (определение~\ref{oval}) на случай  произвольной выпуклой фигуры $F$. Следует только поменять пункт {\it г)} на следующий:

\vskip0.2cm

{\it  \v г) нижний радиус кривизны кривой, ограничивающей $F$, монотонно возрастают  от $A$ до $B$.}

\vskip0.2cm

Это свойство можно переформулировать так: если вы едете на машине от $A$ до $B$ по границе $F$, то вам придётся всё время крутить руль только вправо. Ваши колёса всё время будут повёрнуты влево, иначе $F$ была бы не выпуклой, но сам руль вам придётся крутить вправо, уменьшая угол между передними и задними колёсами. 

В случае, когда граница $F$ --- гладкая кривая, условие {\it\v  г)} совпадает с монотонным возрастанием обычного радиуса кривизны от $A$ до $B$.

Мы, наконец, готовы сформулировать главный результат этой статьи:

\begin{thm}{Теорема}\label{thm-4}
Выпуклая фигура  уникальносоставлена тогда и только тогда, когда она овальна.
\end{thm}

В доказательстве нам потребуется понятие «профиль выпуклой фигуры».
Пусть $F$ есть выпуклая фигура. 
Обозначим через $\rho_F(r)$ угловую меру участка границы $F$ с нижним радиусом кривизны $\le r$.
Функция $\rho_F:\RR_+\to [0,2{\cdot}\pi]$ называется профилем фигуры $F$.

Например, профиль  круга радиуса $R$ есть
$$\rho_K(r)=\left[\begin{matrix}
 0     & \mbox{если} & r<R\\ 
 2{\cdot}\pi & \mbox{если} & r\ge R,
 \end{matrix}\right.$$
 а профиль любого многоугольника тождественно равен $2{\cdot}\pi$.

В случае круговых многоугольников профиль границы есть ступеньчатая функция. Равенство профилей границ круговых многоугольников равносильно стабильной равносоставленности их границ. 
В общем случае это уже не так. 
Тем не менее, «равенство профилей»  выпуклых фигур работает почти так же, как  «стабильная равносоставленность границ», т.е. верен следующий аналог утверждения 1:

\begin{thm}{Утверждение}
Следующие два свойства пары выпуклых фигур $F$ и $G$ равносильны:
 
 (а)  $F$ и $G$ имеют равные профили границ и равную площадь 
 
 (б) для любого $\eps>0$ найдётся пара фигур $F'$ и $G'$, $\eps$-близкие по Хаусдорфу соответственно к $F$ и $G$ и такие, что $F\sim G'$ и $G\sim F'$. 
\end{thm}

Кроме того, легко видеть, что в  классе выпуклых фигур с равным профилем (не равным тождественно $2{\cdot}\pi$) найдётся едиственный, с точностью до конгруэнтости, овал.

Также как в случае круговых многоугольников, теорема \ref{thm-4} 
следует из следующего варианта леммы 
\ref{lm-4}: 

\begin{thm}{Лемма}
Если выпуклая фигура $F$ имеет профиль, равный профилю овала $V$, то $S(F)\geqslant S(V)$. Более того, если $S(F)=S(V)$ то  $F\cong V$.
\end{thm}

\parit{Набросок доказательства.} Произвольную выпуклую фигуру $F$ можно приблизить последовательностью круговых 
многоугольников $M_i$, так что профили $M_i$ сходятся к профилю $F$. (Это почти равносильно тому, что измеримую функцию на отрезке можно приблизить (по мере) ступенчатой функцией.)

Пусть $V_i$ обозначают овалы, соответствующие $M_i$.
Нетрудно видеть, что $V_i$ сходятся по Хаусдорфу к некоторому овалу $V$. Более того, полученный овал $V$ имеет тот же профиль, что $F$. 
Из леммы \ref{lm-4} получаем 
$S(V_i)\le S(M_i)$
и, перейдя к пределу,
$S(V)\le S(F)$.
Остаётся только показать, что в случае равенства  $V$ и $ F$ конгруэнтны.
Для этого надо слегка уточнить лемму~\ref{lm-4}:

\begin{thm}{Лемма}
Для любого $\eps>0$ существует $\delta>0$ такое, что, если $M$ есть круговой
многоугольник, и $V$ есть соответствующий ему овал, то $d_H(V,M)>\eps{\cdot} h$ влечёт $S(M)-S(V)>\delta{\cdot} h^2$, где $h$ обозначает ширину овала $V$.  
\end{thm}

Доказательство получается довольно простым анализом каждого построения в доказательстве леммы~\ref{lm-4}.\qed

\section{Продвинутые упражнения}

\parbf{Аффинная уникальносоставленность.}
Понятие равносоставленности можно также рассматривать для других групп преобразований плоскости, как, например, для группы параллельных переносов или группы подобий. Довольно интересным случаем является группа эквиафинных преобразований, т.е. афинных преобразований сохраняющих площадь. В этом случае, такими  {аффинно-уникальносоставленными} фигурами являются только эллипсы. Попробуйте это доказать

Для решения этого упражнения вам пригодится знать, что такое «аффинная длина», для этого мы советуем почитать книжку Фейеша Тота \cite{F}. 

\parbf{Уникальносоставленные наборы фигур.} Уникальносоставленность можно обобщить на наборы из $n$ выпуклых фигур: 
{\it набор из $n$ выпуклых фигур называется уникальносоставленным, если из равносоставленности этого набора второму набору из $n$ выпуклых фигур следует, что каждая из фигур первого набора конгруэнтна одной из фигур второго}. 

Докажите, что набор из $n$  выпуклых фигур уникальносоставлен тогда и только тогда, когда он состоит из конгруэнтных между собой овалов.

\parbf{Другое определение равносоставленности.} Можно рассматривать другое определение равносоставленности, разрешив разрезать по произвольным Жордановым дугам (т.е. образам однозначных непрерывных отображений отрезка), но как было доказанно в \cite{DKH}, для выпуклых фигур эти определения равносильны, т.е. если две выпуклые фигуры «равносоставлены» по этому новому определению, то они равносоставлены и в нашем смысле. 

Это не так просто, как может показаться; попробуйте доказать это сами.

\parbf{Уникальносоставленные тела.}

\begin{itemize}
\item[а)]
Докажите, что для любого выпуклого тела в $\RR^3$ найдётся произвольно близкое уникальносоставленное выпуклое тело.
\item[б)]
Попробуйте найти точный смысл слов «почти все» и доказать, что «почти все выпуклые тела в трёхмерном пространстве являются уникальносоставленными». Для этого вам придётся узнать, что такое категория Бэра; это можно сделать почитав книжку \cite{O}.
\end{itemize}

\parbf{Равносоставленность неограниченных фигур.} Эту задачу предложил нам А.~В.~Гиль, другой участник уже упомянутого математического кружка. Мы её решили и получили массу удовольствия, теперь предлагаем вам: 

Равносоставленность можно определить для неограниченных фигур, допуская разрезания вдоль конечного числа лучей и отрезков прямых. 

\begin{itemize}
\item[а)] Сформулируйте и докажите аналог теоремы Бойаи --- Гервина для {\it бесконечных многоугольников}, т.е. для фигур ограниченных конечным числом лучей и отрезков. Иначе говоря, найдите необходимое и достаточное условие равносоставленности двух бесконечных многоугольников.

\item[б)] Докажите, что если бесконечная выпуклая фигура $F$ уникальносоставлена, то из каждой её точки исходит ровно один луч целиком содержащийся в $F$.

\item[в)] Опишите все бесконечные уникальносоставленные выпуклые фигуры.
\end{itemize}

\end{document}